\documentclass[11pt]{article}
\usepackage{amsmath}
\usepackage{amssymb}

\usepackage{theorem}
\numberwithin{equation}{section}

\newtheorem{theorem}{Theorem}[section]

\newtheorem{lemma}[theorem]{Lemma}

\begin{document}
\title{Global well-posedness for the defocusing, quintic nonlinear Schr{\"o}dinger equation in one dimension}
\date{\today}
\author{Ben Dodson}
\maketitle

\noindent \textbf{Abstract:} In this paper, we prove global well-posedness for low regularity data for the one dimensional quintic defocusing nonlinear Schr{\"o}dinger equation. We show that a unique solution exists for $u_{0} \in H^{s}(\mathbf{R})$, $s > \frac{8}{29}$. This improves the result in \cite{DPST1}, which proved global well-posedness for $s > \frac{1}{3}$. The main new argument is that we obtain almost Morawetz estimates with improved error.

\section{Introduction}
In this paper we study the initial value problem for the quintic, defocusing, nonlinear Schr{\"o}dinger equation in one dimension,

\begin{equation}\label{1.1}
\aligned
iu_{t} + \Delta u &= |u|^{4} u, \\
u(0,x) &= u_{0} \in H^{s}(\mathbf{R}).
\endaligned
\end{equation}

\noindent This is an $L^{2}$-critical equation. By the results of \cite{CaWe1}, this equation has a local solution on some $[0, T]$, $T(\| u_{0} \|_{H^{s}(\mathbf{R})}) > 0$, when $s > 0$. If a solution to $(\ref{1.1})$ fails to be global and only exists on $[0, T_{\ast})$, $T_{\ast} < \infty$, then

\begin{equation}\label{1.1.1}
\lim_{t \nearrow T_{\ast}} \| u(t) \|_{H^{s}(\mathbf{R})} = \infty.
\end{equation}

\noindent $\cite{CaWe}$ proved

\begin{equation}\label{1.2}
M(u(t)) = \int |u(t,x)|^{2} dx = M(u(0)),
\end{equation}

\begin{equation}\label{1.3}
E(u(t)) = \frac{1}{2} \int |\nabla u(t,x)|^{2} dx + \frac{1}{6} \int |u(t,x)|^{6} dx = E(u(0)),
\end{equation}

\noindent are conserved, giving global well-posedness for $u_{0} \in H^{1}(\mathbf{R})$. The regularity necessary for global well-posedness has since been loweredto $s > 1/3$, (see $\cite{DPST1}$). In this paper we will prove

\begin{theorem}\label{t1.1}
$(\ref{1.1})$ is globally well-posed for all $u_{0} \in H^{s}(\mathbf{R})$, $s > \frac{8}{29}$. Moreover,

\begin{equation}\label{1.4}
\sup_{t \in [0, T]} \| u(t) \|_{H^{s}(\mathbf{R})} \lesssim (1 + T)^{\frac{4(1 - s)s}{29s - 8}+}.
\end{equation}
\end{theorem}

\noindent $\cite{DPST1}$ used the I-method, a method that we will utilize in this paper as well. The I-method was first introduced for the defocusing, cubic initial value problem (see $\cite{CKSTT1}$).\vspace{5mm}

\noindent In $\S 2$, we will start with some preliminary information, including the Strichartz estimates, Littlewood-Paley theory, a description of the I-method, and a local well-posedness result. In $\S 3$, an energy increment will be obtained. In $\S 4$, the almost Morawetz estimates will be proved. In $\S 5$, we will prove the theorem.

\section{Local Well-posedness}
\noindent The proof of local well-posedness makes use of the Strichartz estimates.

\begin{theorem}\label{t2.0.1}
\noindent A pair $(p,q)$ is called an admissible pair if $\frac{2}{p} + \frac{1}{q} = \frac{1}{2}$. If $(p,q)$ and $(\tilde{p}, \tilde{q})$ are admissible pairs and and $u(t,x)$ solves

\begin{equation}\label{2.0.1}
\aligned
i u_{t} + \Delta u &= F(t), \\
u(0,x) &= u_{0},
\endaligned
\end{equation}

\noindent then

\begin{equation}\label{2.0.1}
\| u(t,x) \|_{L_{t}^{p} L_{x}^{q}(J \times \mathbf{R})} \lesssim \| u_{0} \|_{L^{2}(\mathbf{R})} + \| F(t) \|_{L_{t}^{\tilde{p}'} L_{x}^{\tilde{q}'}(J \times \mathbf{R})}.
\end{equation}
\end{theorem}

\noindent \emph{Proof:} See \cite{Tao}. $p'$ denotes the Lebesgue exponent $\frac{p}{p - 1}$.\vspace{5mm}

\noindent The Strichartz space will be defined by the norm

\begin{equation}\label{2.0.1.0}
\| u \|_{S^{0}(J \times \mathbf{R})} = \sup_{(p,q) \text{ admissible}} \| u \|_{L_{t}^{p} L_{x}^{q}(J \times \mathbf{R})}.
\end{equation}

\noindent The space $N^{0}(J \times \mathbf{R})$ is the dual space to $S^{0}(J \times \mathbf{R})$. See \cite{Tao} for more details.\vspace{5mm}

\noindent We will also make use of the Littlewood-Paley decomposition. Suppose $\phi(x)$ is a smooth function,

\begin{equation}\label{2.0.1.1}
\phi(x) = \left\{
            \begin{array}{ll}
              1, & \hbox{$|x| \leq 1/2$;} \\
              0, & \hbox{$|x| > 1$.}
            \end{array}
          \right.
\end{equation}

\begin{equation}\label{2.0.2}
\aligned
\mathcal F(P_{\leq N} u) = \phi(\frac{\xi}{N}) \hat{u}(\xi), \\
\mathcal F(P_{> N} u) = (1 - \phi(\frac{\xi}{N})) \hat{u}(\xi), \\
\mathcal F(P_{N} u) = P_{\leq N} u - P_{\leq \frac{N}{2}} u.
\endaligned
\end{equation}

\noindent For convenience, let $u_{N} = P_{N} u$, similarly for $u_{\leq N}$, $u_{> N}$.\vspace{5mm}

\noindent The I - operator is a Fourier multiplier,

\begin{equation}\label{2.0.3}
I_{N} : H^{s}(\mathbf{R}) \rightarrow H^{1}(\mathbf{R}),
\end{equation}

\begin{equation}\label{2.0.4}
\widehat{I_{N} f}(\xi) = m_{N}(\xi) \hat{f}(\xi),
\end{equation}

\begin{equation}\label{2.0.5}
m_{N}(\xi) = \left\{
               \begin{array}{ll}
                 1, & \hbox{$|\xi| \leq N$;} \\
                 (\frac{N}{|\xi|})^{s - 1}, & \hbox{$|\xi| > 2N$.}
               \end{array}
             \right.
\end{equation}

\begin{equation}\label{2.0.6}
\aligned
\| Iu \|_{H^{1}(\mathbf{R})} &\lesssim N^{1 - s} \| u \|_{H^{s}(\mathbf{R})},\\
\| u \|_{H^{s}(\mathbf{R})} &\lesssim \| Iu \|_{H^{1}(\mathbf{R})},
\endaligned
\end{equation}

\noindent therefore, controlling $E(Iu(t))$ gives control of $\| u(t) \|_{H^{s}(\mathbf{R})}$. For the rest of the paper, $If$ denotes $I_{N} f$, and the presence of an $N$ is implied.

\begin{lemma}\label{l2.0}
Let I be a compact time interval, $t_{0} \in I$, $N > 0$, and suppose $u_{1}$, $u_{2}$ are two solutions to $(\ref{2.0.1})$ such that $u_{j}(t)$ has Fourier support in the region $\{ |\xi_{j}| \leq N \}$ for $j = 1, 2$. Suppose also that the Fourier supports of $u_{1}, u_{2}$ are separated by at least $\geq cN$. Then for any $q > 2$,

\begin{equation}\label{3.4.1}
\| u_{1} u_{2} \|_{L_{t,x}^{q}(I \times \mathbf{R})} \lesssim N^{1 - 3/q} \| u_{1} \|_{S_{\ast}^{0}(I \times \mathbf{R})} \| u_{2} \|_{S_{\ast}^{0}(I \times \mathbf{R})},
\end{equation}

\noindent where

\begin{equation}\label{3.4.2}
\| u \|_{S_{\ast}^{0}(I \times \mathbf{R})} = \| u_{0} \|_{L^{2}(\mathbf{R})} + \| (i\partial_{t} + \Delta) u \|_{L_{t,x}^{6/5}(I \times \mathbf{R})}.
\end{equation}
\end{lemma}

\noindent \emph{Proof:} See \cite{TaoT}.\vspace{5mm}

\noindent To this end, let $q = 2 + \delta$, $\frac{1}{p} = \frac{1}{2} - \frac{1}{2 + \delta}$.\vspace{5mm}

\noindent In proving theorem $\ref{t1.1}$, we will make use of a linear-nonlinear decomposition. See \cite{R} for the linear-nonlinear decomposition for the defocusing, semilinear wave equation, \cite{D1} for the linear-nonlinear decomposition used for the three dimensional cubic defocusing nonlinear Schr{\"o}dinger equation.\vspace{5mm}

\begin{theorem}\label{t2.1}
If

\begin{equation}\label{2.1}
\| \nabla Iu_{0} \|_{L^{2}(\mathbf{R})} \leq 1
\end{equation}

\noindent and for some $\epsilon > 0$ sufficiently small,

\begin{equation}\label{2.2}
\| Iu \|_{L_{t}^{16/3} L_{x}^{8}(J \times \mathbf{R})} \leq \epsilon,
\end{equation}

\noindent then

\begin{equation}\label{2.3}
\| \langle \nabla \rangle Iu \|_{S^{0}(J \times \mathbf{R})} \lesssim 1.
\end{equation}

\noindent Moreover, the solution has the form

\begin{equation}\label{2.3.1}
e^{it \Delta} u_{0} + u^{nl}(t),
\end{equation}

\begin{equation}\label{2.3.2}
\| P_{> cN} \langle \nabla \rangle I u^{nl} \|_{S^{0}(J \times \mathbf{R})} \lesssim \frac{1}{N^{1/2-}}.
\end{equation}
\end{theorem}

\noindent \emph{Proof:} The solution obeys the Duhamel formula,

\begin{equation}\label{2.4}
Iu(t,x) = e^{it \Delta} u_{0} + \int_{0}^{t} e^{i(t - \tau) \Delta} I(|u(\tau)|^{4} u(\tau)) d\tau.
\end{equation}

\noindent By the Strichartz estimates,

\begin{equation}\label{2.5}
\| (1 - I) u \|_{L_{t}^{16/3} L_{x}^{8}(J \times \mathbf{R})} \lesssim \frac{1}{N} \| \langle \nabla \rangle Iu \|_{S^{0}(J \times \mathbf{R})},
\end{equation}

\begin{equation}\label{2.5.1}
\| \langle \nabla \rangle Iu \|_{S^{0}(J \times \mathbf{R})} \lesssim \| \langle \nabla \rangle Iu_{0} \|_{L^{2}(\mathbf{R})} + \| \langle \nabla \rangle Iu \|_{S^{0}(J \times \mathbf{R})} \| u \|_{L_{t}^{16/3} L_{x}^{8}(J \times \mathbf{R})}^{4}
\end{equation}

\begin{equation}\label{2.6}
\lesssim \| \langle \nabla \rangle Iu \|_{S^{0}(J \times \mathbf{R})} \lesssim \| \langle \nabla \rangle Iu_{0} \|_{L^{2}(\mathbf{R})} + \| \langle \nabla \rangle Iu \|_{S^{0}(J \times \mathbf{R})} (\epsilon^{4} + \frac{\| \langle \nabla \rangle Iu \|_{S^{0}(J \times \mathbf{R})}^{4}}{N^{4}}).
\end{equation}

\noindent Therefore, by the continuity method,

\begin{equation}\label{2.7}
\| \langle \nabla \rangle Iu \|_{S^{0}(J \times \mathbf{R})} \lesssim 1.
\end{equation}

\noindent This takes care of $(\ref{2.3})$. Next, we remark that this also proves

\begin{equation}\label{2.7.1}
\| \langle \nabla \rangle I(|u|^{4} u) \|_{L_{t,x}^{6/5}(J \times \mathbf{R})} \lesssim \| \langle \nabla \rangle Iu \|_{S^{0}(J \times \mathbf{R})}^{5} \lesssim 1.
\end{equation}

\noindent To estimate the nonlinearity,

\begin{equation}\label{2.8}
\nabla Iu^{nl}(t) = \int_{0}^{t} \nabla e^{i(t - \tau) \Delta} I(|u(\tau)|^{4} u(\tau)) d\tau.
\end{equation}

\begin{equation}\label{2.9}
\aligned
&I(|u(\tau)|^{4} u(\tau)) = I(|u_{\leq \frac{cN}{10}}(\tau)|^{4} u_{\leq \frac{cN}{10}}(\tau)) + I(O((u_{\leq \frac{cN}{100}}(\tau))^{4} (u_{> \frac{cN}{10}}(\tau)))) \\ &+ I(O((u_{> \frac{cN}{10}}(\tau))(u_{> \frac{cN}{100}}(\tau)) u(\tau)^{3})).
\endaligned
\end{equation}

\noindent The first term, $I(|u_{\leq \frac{cN}{10}}(\tau)|^{4} u_{\leq \frac{cN}{10}}(\tau))$ is supported on $|\xi| \leq \frac{cN}{2}$. To estimate the second term, use the bilinear estimates,

$$\| \nabla I(O(u_{> \frac{cN}{10}} (u_{\leq \frac{cN}{100}})^{4})) \|_{L_{t}^{1} L_{x}^{2}(J \times \mathbf{R})}$$ $$\lesssim \| (\nabla Iu_{> \frac{cN}{10}})(u_{\leq \frac{cN}{100}}) \|_{L_{t,x}^{q}} \| u_{\leq \frac{cN}{100}} \|_{L_{t,x}^{p}} \| u_{\leq \frac{cN}{100}} \|_{L_{t}^{4} L_{x}^{\infty}}^{2} \lesssim \frac{1}{N^{1/2-}}.$$

\noindent Finally,

$$\| \nabla I(O((u_{> \frac{cN}{10}}(\tau))(u_{> \frac{cN}{100}}(\tau)) u(\tau)^{3})) \|_{L_{t}^{1} L_{x}^{2}(J \times \mathbf{R})}$$ $$\lesssim \| \nabla Iu \|_{L_{t}^{4} L_{x}^{\infty}(J \times \mathbf{R})} \| u_{> \frac{cN}{100}} \|_{L_{t}^{16/3} L_{x}^{8}(J \times \mathbf{R})} \| u \|_{L_{t}^{16/3} L_{x}^{8}(J \times \mathbf{R})}^{3}$$

$$\lesssim \frac{\epsilon^{3}}{N} \| \langle \nabla \rangle Iu \|_{S^{0}(J \times \mathbf{R})}^{2} \lesssim \frac{1}{N}.$$

\noindent The last estimate follows from

\begin{equation}\label{2.10}
\aligned
\| u_{> \frac{cN}{100}} \|_{L_{t}^{16/3} L_{x}^{8}(J \times \mathbf{R})} \lesssim \sum_{\frac{cN}{100} \leq N_{j}} \| P_{N_{j}} u \|_{L_{t}^{16/3} L_{x}^{8}(J \times \mathbf{R})} \\
\lesssim \| \langle \nabla \rangle Iu \|_{S^{0}(J \times \mathbf{R})} \sum_{\frac{cN}{100} \leq N_{j}} \frac{1}{N_{j}^{s} N^{1 - s}} \lesssim \frac{1}{N} \| \langle \nabla \rangle Iu \|_{S^{0}(J \times \mathbf{R})}.
\endaligned
\end{equation}

\noindent $\Box$

\section{Energy Increment}
\noindent In this section we prove almost conservation of the modified energy $E(Iu(t))$.

\begin{theorem}\label{t3.1}
Let

\begin{equation}\label{3.1}
E(Iu(t)) = \frac{1}{2} \int |\nabla Iu(t,x)|^{2} dx + \frac{1}{6} \int |Iu(t,x)|^{6} dx.
\end{equation}

\noindent If $J$ is an interval where a solution $u(t,x)$ of $(\ref{1.1})$ exists, $\| u \|_{L_{t}^{16/3} L_{x}^{8}(J \times \mathbf{R})} \leq \epsilon$, $E(Iu_{0}) \leq 1$, then

\begin{equation}\label{3.1.1}
\sup_{t_{1}, t_{2} \in J} |E(Iu(t_{1})) - E(Iu(t_{2}))| \lesssim \frac{1}{N^{3/2-}} \| P_{> cN} \nabla Iu \|_{L_{t}^{4} L_{x}^{\infty}(J \times \mathbf{R})}^{2} + \frac{1}{N^{2-}},
\end{equation}

\noindent where $c > 0$ is some small constant.
\end{theorem}

\noindent \emph{Proof:}

\begin{equation}\label{3.2}
\frac{d}{dt} E(Iu(t)) = Re \int (\overline{Iu_{t}(t,x)}) [I(|u(t,x)|^{4} u(t,x)) - |Iu(t,x)|^{4} Iu(t,x)] dx.
\end{equation}

\noindent Taking the Fourier transform, let $\Sigma = \{ \xi_{1} + ... + \xi_{6} = 0 \}$, $d \xi$ is the Lebesgue measure on the hyperplane, using the fact that

\begin{equation}\label{3.2.1}
Iu_{t} = iI \Delta u - iI(|u|^{4} u),
\end{equation}

\begin{equation}\label{3.3}
\aligned
\frac{d}{dt} E(Iu(t)) = -Re \int_{\Sigma} (i |\xi_{1}|^{2} \widehat{\overline{Iu}}(t,\xi_{1})) [1 - \frac{m(\xi_{2} + ... + \xi_{6})}{m(\xi_{2}) m(\xi_{3}) \cdots m(\xi_{6})}] \\ \times \widehat{Iu}(t,\xi_{2}) \widehat{\overline{Iu}}(t,\xi_{3}) \widehat{Iu}(t,\xi_{4}) \widehat{\overline{Iu}}(t,\xi_{5}) \widehat{Iu}(t,\xi_{6}) d\xi
\endaligned
\end{equation}

\begin{equation}\label{3.4}
\aligned
- Re \int_{\Sigma} (i  \widehat{\overline{I(|u|^{4} u)}}(t,\xi_{1})) [1 - \frac{m(\xi_{2} + ... + \xi_{6})}{m(\xi_{2}) m(\xi_{3}) \cdots m(\xi_{6})}] \\
 \times \widehat{Iu}(t,\xi_{2}) \widehat{\overline{Iu}}(t,\xi_{3}) \widehat{Iu}(t,\xi_{4}) \widehat{\overline{Iu}}(t,\xi_{5}) \widehat{Iu}(t,\xi_{6}) d\xi.
\endaligned
\end{equation}

\noindent We will estimate $(\ref{3.3})$ and $(\ref{3.4})$ separately by making a Littlewood-Paley decomposition and consider several cases separately. Without loss of generality let $N_{2} \geq N_{3} \geq N_{4} \geq N_{5} \geq N_{6}$.\vspace{5mm}

\noindent \textbf{The term $(\ref{3.3})$:}\vspace{5mm}

\noindent When estimating this term, we will frequently use the bilinear estimate $(\ref{3.4.1})$.

\noindent \textbf{Case 1, $N_{2} << N$:} In this case, $m(\xi_{i}) \equiv 1$, so

\begin{equation}\label{3.5}
1 - \frac{m(\xi_{2} + ... + \xi_{6})}{m(\xi_{2}) \cdots m(\xi_{6})} \equiv 0.
\end{equation}

\noindent \textbf{Case 2, $N_{2} \gtrsim N >> N_{3}$:} By the fundamental theorem of calculus,

\begin{equation}\label{3.6}
|1 - \frac{m(\xi_{2} + ... + \xi_{6})}{m(\xi_{2}) \cdots m(\xi_{6})}| \lesssim \frac{N_{3}}{N_{2}}.
\end{equation}

\noindent Recall that $q = 2 + \delta$, $\frac{1}{p} = \frac{1}{2} - \frac{1}{q}$. $N_{1} \sim N_{2}$, so $$\sum_{N \lesssim N_{1} \sim N_{2}} \sum_{N_{6} \leq N_{5} \leq N_{4} \leq N_{3} << N} \frac{N_{1}}{N_{2}^{2}} \| (P_{N_{1}} \nabla Iu)(P_{N_{3}} Iu) \|_{L_{t,x}^{q}}$$

$$\times \| (P_{N_{2}} \nabla Iu)(P_{N_{4}} Iu) \|_{L_{t,x}^{q}} \|P_{N_{5}} Iu \|_{L_{t,x}^{p}} \| P_{N_{6}} Iu \|_{L_{t,x}^{p}}$$

$$\lesssim \sum_{N \lesssim N_{1} \sim N_{2}} \sum_{N_{6} \leq N_{5} \leq N_{4} \leq N_{3} << N} \frac{N_{1}^{1/2+}}{N_{2}^{5/2-}} \frac{1}{\langle N_{3} \rangle} \frac{1}{\langle N_{4} \rangle \langle N_{5} \rangle^{1/2-} \langle N_{6} \rangle^{1/2-}} \lesssim \frac{1}{N^{2-}}.$$

\noindent \textbf{Case 3, $N_{2} \geq N_{3} \gtrsim N >> N_{4}$:} In this case, estimate the multiplier by

\begin{equation}\label{3.7}
|1 - \frac{m(\xi_{1})}{m(\xi_{2}) m(\xi_{3})}| \lesssim \frac{m(\xi_{1})}{m(\xi_{2}) m(\xi_{3})}.
\end{equation}

\noindent Consider three subcases separately.\vspace{5mm}

\noindent \emph{Case 3(a), $N_{3} << N_{1} \sim N_{2} $:} $$\sum_{N \lesssim N_{3} \lesssim N_{1} \sim N_{2}} \sum_{N_{6} \leq N_{5} \leq N_{4} << N} \frac{1}{m(N_{3})} \frac{N_{1}}{N_{2}} \| (P_{N_{1}} \nabla Iu) (P_{N_{3}} Iu) \|_{L_{t,x}^{q}}$$

$$\| (P_{N_{2}} \nabla Iu) (P_{N_{4}} Iu) \|_{L_{t,x}^{q}} \| P_{N_{5}} Iu \|_{L_{t,x}^{p}} \| P_{N_{6}} Iu \|_{L_{t,x}^{p}}$$

$$\lesssim \sum_{N \lesssim N_{3} << N_{1} \sim N_{2}} \sum_{N_{6} \leq N_{5} \leq N_{4} << N}  \frac{N_{1}^{1/2}}{N_{2}^{3/2-}} \frac{1}{N_{3}^{s} N^{1 - s}} \frac{1}{\langle N_{4} \rangle \langle N_{5} \rangle^{1/2-} \langle N_{6} \rangle^{1/2-}} \lesssim \frac{1}{N^{2-}}. $$\vspace{5mm}

\noindent \emph{Case 3(b), $N_{1} << N_{2} \sim N_{3}$:}

$$\sum_{N \lesssim N_{2} \sim N_{3}} \sum_{N_{6} \leq N_{5} \leq N_{4} << N, N_{1} << N_{2}} \frac{m(N_{1})}{m(N_{2}) m(N_{3})} \frac{N_{1}}{N_{2}} \| (P_{N_{2}} \nabla Iu) (P_{N_{1}} \nabla Iu) \|_{L_{t,x}^{q}}$$

$$\| (P_{N_{3}} Iu) (P_{N_{4}} Iu) \|_{L_{t,x}^{q}} \| P_{N_{5}} Iu \|_{L_{t,x}^{p}} \| P_{N_{6}} Iu \|_{L_{t,x}^{p}}$$

$$\lesssim \sum_{N \lesssim N_{3} \lesssim N_{1} \sim N_{2}} \sum_{N_{6} \leq N_{5} \leq N_{4} << N} \frac{1}{N_{1}^{1/2-} N_{2}^{1/2-}} \frac{1}{N_{2}^{s} N^{1 - s}} \frac{N_{1}^{s} N^{1 - s}}{N_{3}^{s} N^{1 - s}}$$ $$\times \frac{1}{\langle N_{4} \rangle \langle N_{5} \rangle^{1/2-} \langle N_{6} \rangle^{1/2-}} \lesssim \frac{1}{N^{2-}}.$$\vspace{5mm}

\noindent \emph{Case 3(c), $N_{1} \sim N_{2} \sim N_{3}$:} In this case,

$$\frac{m(N_{1})}{m(N_{2}) m(N_{3})} \sim \frac{1}{m(N_{3})}.$$

$$\sum_{N \lesssim N_{1} \sim N_{2} \sim N_{3}} \sum_{N_{6} \leq N_{5} \leq N_{4} << N} \frac{N_{1}}{N_{2}} \| (P_{N_{1}} \nabla Iu)(P_{N_{4}} Iu) \|_{L_{t,x}^{q}}$$

$$\| P_{N_{2}} \nabla Iu \|_{L_{t}^{4} L_{x}^{\infty}} \| P_{N_{3}} Iu \|_{L_{t}^{4} L_{x}^{\infty}} \| P_{N_{5}} Iu \|_{L_{t}^{\infty} L_{x}^{2}} \| P_{N_{6}} Iu \|_{L_{t,x}^{p}}$$

$$\lesssim \sum_{N \lesssim N_{2} \sim N_{3}} \frac{1}{N_{2}^{1/2-} N^{1 - s} N_{3}^{s}} \| P_{N_{2}} \nabla Iu \|_{L_{t}^{4} L_{x}^{\infty}} \| P_{N_{3}} \nabla Iu \|_{L_{t}^{4} L_{x}^{\infty}}$$ $$\times \sum_{N_{6} \leq N_{5} \leq N_{4} << N} \frac{1}{\langle N_{4} \rangle \langle N_{5} \rangle \langle N_{6} \rangle^{1/2-}} \lesssim \frac{1}{N^{3/2-}} \| P_{> cN} Iu \|_{L_{t}^{4} L_{x}^{\infty}(J \times \mathbf{R})}^{2} .$$\vspace{5mm}

\noindent \textbf{Case 4, $N_{2} \geq N_{3} \geq N_{4} \gtrsim N$:}\vspace{5mm}

\noindent \emph{Case 4(a), $N_{1} \sim N_{2}$:} In this case $$|1 - \frac{m(\xi_{1})}{m(\xi_{2}) \cdots m(\xi_{6})} | \lesssim \frac{1}{m(\xi_{3}) \cdots m(\xi_{6})}.$$

$$\sum_{N \lesssim N_{1} \sim N_{2}} \sum_{N_{6} \leq N_{5} \leq N_{4} \leq N_{3}} \frac{N_{1}}{N_{2}} \| P_{N_{1}} \nabla Iu \|_{L_{t,x}^{6}} \| P_{N_{2}} \nabla Iu \|_{L_{t,x}^{6}}$$ $$\times \| P_{N_{3}} u \|_{L_{t,x}^{6}} \| P_{N_{4}} u \|_{L_{t,x}^{6}} \| P_{N_{5}} u \|_{L_{t,x}^{6}} \| P_{N_{6}} u \|_{L_{t,x}^{6}}$$

$$\lesssim \sum_{N \lesssim N_{1} \sim N_{2}} \sum_{N \lesssim N_{4} \leq N_{3}} \sum_{N_{6} \leq N_{5} \leq N} \frac{N_{1}}{N_{2}} \| P_{N_{1}} \nabla Iu \|_{L_{t,x}^{6}} \| P_{N_{2}} \nabla Iu \|_{L_{t,x}^{6}}$$

$$\times \sum_{N \lesssim N_{4} \leq N_{3}} \frac{1}{N_{3}^{s} N^{1 - s}} \frac{1}{N_{4}^{s} N^{1 - s}} \frac{1}{m(N_{5}) N_{5}} \frac{1}{m(N_{6}) N_{6}} \lesssim \frac{1}{N^{2-}}.$$\vspace{5mm}

\noindent \emph{Case 2, $N_{1} << N_{2}$:} In this case $N_{2} \sim N_{3}$. $$\sum_{N \lesssim N_{2} \sim N_{3}} \sum_{N \lesssim N_{4}; N_{6} \leq N_{5} \leq N_{4}; N_{1} << N_{2}} \frac{N_{1}}{N_{2}} \frac{m(N_{1})}{m(N_{2})} \| P_{N_{1}} \nabla Iu \|_{L_{t,x}^{6}} \| P_{N_{2}} \nabla Iu \|_{L_{t,x}^{6}}$$

$$\| P_{N_{3}} u \|_{L_{t,x}^{6}} \| P_{N_{4}} u \|_{L_{t,x}^{6}} \| P_{N_{5}} u \|_{L_{t,x}^{6}} \| P_{N_{6}} u \|_{L_{t,x}^{6}}$$
$$\lesssim \sum_{N \lesssim N_{2} \sim N_{3}} \sum_{N \lesssim N_{4}; N_{6} \leq N_{5} \leq N_{4}; N_{1} << N_{2}} \frac{N_{1}}{N_{2}} \frac{m(N_{1})}{m(N_{2})} \frac{1}{N_{3}^{s} N^{1 - s}} \frac{1}{N_{4}^{s} N^{1 - s}}$$ $$\times \frac{1}{N_{5} m(N_{5})} \frac{1}{N_{6} m(N_{6})}  \lesssim \frac{1}{N^{2-}}.$$

\noindent This takes care of $(\ref{3.3})$.\vspace{5mm}

\noindent \textbf{The term $(\ref{3.4})$:}  Recall that this is the $10$-linear term

\begin{equation}\label{3.8}
Re \int_{0}^{t} \int \overline{i I(|u|^{4} u)} [I(|u|^{4} u) - |Iu|^{4} (Iu)] dx dt.
\end{equation}

\noindent The term $I(|u|^{4} u)$ poses a slight technical problem. Ideally, this term would be placed in $L_{t,x}^{6}$, and we would then repeat the analysis used in $(\ref{3.3})$. However, in general this is not possible, so instead let $u_{l} = P_{\leq \frac{N}{30}} u$, $u_{h} = P_{> \frac{N}{30}} u$,

\begin{equation}\label{3.9}
F(t,x) = \overline{I(|u|^{4} u)} - I(\overline{O(u_{h}^{4} u_{l})}) - I(\overline{|u_{h}|^{4} u_{h}}),
\end{equation}

\noindent where $O(u_{h}^{4} u_{l})$ consists of those terms in $|u_{l} + u_{h}|^{4} (u_{l} + u_{h})$ consisting of four $u_{h}$ terms and one $u_{l}$ term.

\begin{equation}\label{3.10}
\aligned
\| F(t,x) \|_{L_{t}^{4} L_{x}^{\infty}(J \times \mathbf{R})} \lesssim \| \langle \nabla \rangle F(t,x) \|_{L_{t}^{4} L_{x}^{1}(J \times \mathbf{R})} \\ \lesssim \| \langle \nabla \rangle Iu \|_{L_{t}^{4} L_{x}^{\infty}(J \times \mathbf{R})} \| u_{l} \|_{L_{t,x}^{\infty}(J \times \mathbf{R})}^{2} \| u \|_{L_{t}^{\infty} L_{x}^{2}(J \times \mathbf{R})}^{2} \lesssim \| \langle \nabla \rangle Iu \|_{S^{0}(J \times \mathbf{R})}^{5}.
\endaligned
\end{equation}

\noindent Now evaluate

\begin{equation}\label{3.11}
\int_{0}^{t} \int_{\Sigma} \hat{F}(t,\xi_{1}) [1 - \frac{m(\xi_{2} + ... + \xi_{6})}{m(\xi_{2}) \cdots m(\xi_{6})}] \widehat{Iu}(t,\xi_{2}) \widehat{\overline{Iu}}(t,\xi_{3}) \widehat{Iu}(t,\xi_{4}) \widehat{\overline{Iu}}(t,\xi_{5}) \widehat{Iu}(t,\xi_{6}),
\end{equation}

\noindent via a Littlewood - Paley partition of unity and considering several subcases separately. Without loss of generality, let $N_{2} \geq N_{3} \geq ... \geq N_{6}$.\vspace{5mm}

\noindent \textbf{Case 1, $N_{2} << N$:} Once again, $m(\xi_{i}) \equiv 1$, so

\begin{equation}\label{3.12}
1 - \frac{m(\xi_{2} + ... + \xi_{6})}{m(\xi_{2}) \cdots m(\xi_{6})} \equiv 0.
\end{equation}

\noindent \textbf{Case 2, $N_{2} \gtrsim N >> N_{3}$:} In this case, apply the fundamental theorem of calculus,

\begin{equation}\label{3.13}
|1 - \frac{m(\xi_{2} + ... + \xi_{6})}{m(\xi_{2}) \cdots m(\xi_{6})}| \lesssim \frac{N_{3}}{N_{2}}.
\end{equation}

$$\sum_{N \lesssim N_{1} \sim N_{2}} \sum_{N_{6} \leq N_{5} \leq N_{4} \leq N_{3} << N} \frac{N_{3}}{N_{2}} \| P_{N_{1}} F \|_{L_{t}^{4} L_{x}^{\infty}(J \times \mathbf{R})} \| P_{N_{2}} Iu \|_{L_{t}^{4}L_{x}^{\infty}(J \times \mathbf{R})}$$ $$\times \| P_{N_{3}} Iu \|_{L_{t}^{4}L_{x}^{\infty}(J \times \mathbf{R})} \| P_{N_{4}} Iu \|_{L_{t}^{4}L_{x}^{\infty}(J \times \mathbf{R})} \| P_{N_{5}} Iu \|_{L_{t}^{\infty} L_{x}^{2}(J \times \mathbf{R})} \| P_{N_{6}} Iu \|_{L_{t}^{\infty} L_{x}^{2}(J \times \mathbf{R})}$$

$$\lesssim  \sum_{N \lesssim N_{1} \sim N_{2}} \sum_{N_{6} \leq N_{5} \leq N_{4} \leq N_{3} << N} \frac{1}{N_{2}^{2}} \frac{1}{\langle N_{4} \rangle} \frac{1}{\langle N_{5} \rangle} \frac{1}{\langle N_{6} \rangle} \lesssim \frac{1}{N^{2-}}.$$

\noindent \textbf{Case 3, $N_{2} \geq N_{3} \gtrsim N$} In this case make the crude estimate

\begin{equation}\label{3.14}
|1 - \frac{m(\xi_{2} + ... + \xi_{6})}{m(\xi_{2}) \cdots m(\xi_{6})}| \lesssim \frac{1}{m(\xi_{2}) \cdots m(\xi_{6})}.
\end{equation}

$$\sum_{N_{1} \lesssim N_{2}, N \lesssim N_{3} \leq N_{2}} \sum_{N_{6} \leq N_{5} \leq N_{4} \leq N_{3}} \| P_{N_{1}} F \|_{L_{t}^{4} L_{x}^{\infty}(J \times \mathbf{R})} \| P_{N_{2}} u \|_{L_{t}^{4} L_{x}^{\infty}(J \times \mathbf{R})} \| P_{N_{3}} u \|_{L_{t}^{4} L_{x}^{\infty}(J \times \mathbf{R})}$$ $$\times \| P_{N_{4}} u \|_{L_{t}^{4} L_{x}^{\infty} (J \times \mathbf{R})} \| P_{N_{5}} u \|_{L_{t}^{\infty} L_{x}^{2}(J \times \mathbf{R})} \| P_{N_{6}} u \|_{L_{t}^{\infty} L_{x}^{2}(J \times \mathbf{R})}$$

$$\lesssim  \sum_{N_{1} \lesssim N_{2}, N \lesssim N_{3} \leq N_{2}} \sum_{N_{6} \leq N_{5} \leq N_{4} \leq N_{3}} \frac{1}{N_{2}^{s} N^{1 - s}} \frac{1}{N_{3}^{s} N^{1 - s}}$$ $$\times \frac{1}{\langle N_{4} \rangle m(N_{4})} \frac{1}{\langle N_{5} \rangle m(N_{5})} \frac{1}{\langle N_{6} \rangle m(N_{6})} \lesssim \frac{1}{N^{2-}}.$$

\noindent It only remains to consider

\begin{equation}\label{3.15}
\aligned
Re \int_{0}^{t} \int (i \overline{I(|u_{h}|^{4} u_{h})}) [I(|u|^{4} u) - |Iu|^{4} (Iu)] dx dt \\
+ Re \int_{0}^{t} \int i \overline{I(3|u_{h}|^{4} u_{l} + 2 |u_{h}|^{2} u_{h}^{2} \overline{u_{l}})} [I(|u|^{4} u) - |Iu|^{4} (Iu)] dx dt.
\endaligned
\end{equation}

\noindent To evaluate

\begin{equation}\label{3.16}
\aligned
Re \int_{0}^{t} \int (i \overline{I(|u_{h}|^{4} u_{h})}) I(|u|^{4} u) dx dt \\ + Re \int_{0}^{t} \int i \overline{I(3|u_{h}|^{4} u_{l} + 2 |u_{h}|^{2} u_{h}^{2} \overline{u_{l}})} I(|u|^{4} u) dx dt,
\endaligned
\end{equation}

\noindent it is necessary to take advantage of some cancellations.

$$Re \int_{0}^{t} \int (i \overline{I(|u_{h}|^{4} u_{h})}) I(|u_{h}|^{4} u_{h}) dx dt \equiv 0.$$

$$Re \int_{0}^{t} \int i \overline{I(3|u_{h}|^{4} u_{l} + 2 |u_{h}|^{2} u_{h}^{2} \overline{u_{l}})} I(|u_{h}|^{4} u_{h}) dx dt$$ $$+ Re \int_{0}^{t} \int i \overline{I(|u_{h}|^{4} u_{h})} (3|u_{h}|^{4} u_{l} + 2 |u_{h}|^{2} u_{h}^{2} \overline{u_{l}}) dx dt \equiv 0.$$

\noindent It remains to evaluate

\begin{equation}\label{3.17}
\int_{0}^{t} \int I(|u_{h}|^{4} u_{h}) I(O(u_{l}^{2} u^{3})) dx dt + \int_{0}^{t} \int I(O(u_{h}^{4} u_{l})) I(O(u_{l} u^{4})) dx dt.
\end{equation}

\noindent By $(\ref{3.10})$, $$\| I(u_{l}^{2} u^{3}) \|_{L_{t}^{4} L_{x}^{\infty}(J \times \mathbf{R})} \lesssim \| \langle \nabla \rangle Iu \|_{S^{0}(J \times \mathbf{R})}^{5},$$ so

$$\int_{0}^{t} \int I(|u_{h}|^{4} u_{h}) I(O(u_{l}^{2} u^{3})) dx dt \lesssim \| I(u_{l}^{2} u^{3}) \|_{L_{t}^{4} L_{x}^{\infty} (J \times \mathbf{R})} \| u_{h} \|_{L_{t}^{4} L_{x}^{\infty}(J \times \mathbf{R})}^{3} \| u_{h} \|_{L_{t}^{\infty} L_{x}^{2}(J \times \mathbf{R})}^{2}$$ $$\lesssim \frac{1}{N^{5}} \| \langle \nabla \rangle Iu \|_{S^{0}(J \times \mathbf{R})}^{10} \lesssim \frac{1}{N^{5}}.$$

$$\int_{0}^{t} \int I(O(u_{h}^{4} u_{l})) I(O(u_{l} u^{4})) dx dt \lesssim \| \langle \nabla \rangle I(u_{h}^{4} u_{l}) \|_{L_{t}^{2} L_{x}^{1}(J \times \mathbf{R})} \| I(u_{l} u^{4}) \|_{L_{t}^{2} L_{x}^{1}(J \times \mathbf{R})}$$

$$\lesssim \| \langle \nabla \rangle Iu \|_{L_{t}^{\infty} L_{x}^{2}(J \times \mathbf{R})} \| u_{l} \|_{L_{t,x}^{\infty}(J \times \mathbf{R})} \| u_{h} \|_{L_{t,x}^{6}(J \times \mathbf{R})}^{3} \| u_{l} \|_{L_{t,x}^{\infty}(J \times \mathbf{R})} \| u \|_{L_{t,x}^{6}(J \times \mathbf{R})}^{3} \| u \|_{L_{t}^{\infty} L_{x}^{2}(J \times \mathbf{R})}$$

$$\lesssim \frac{1}{N^{3}} \| \langle \nabla \rangle Iu \|_{S^{0}(J \times \mathbf{R})} \lesssim \frac{1}{N^{3}}.$$

\noindent This takes care of $(\ref{3.16})$. To finish estimating $(\ref{3.15})$,

\begin{equation}\label{3.18}
\aligned
Re \int_{0}^{t} \int (i \overline{I(|u_{h}|^{4} u_{h})}) (|Iu|^{4} (Iu) dx dt \\ + Re \int_{0}^{t} \int i \overline{I(3|u_{h}|^{4} u_{l} + 2 |u_{h}|^{2} u_{h}^{2} \overline{u_{l}})} (|Iu|^{4} (Iu) dx dt, \\
\lesssim \| u_{h} \|_{L_{t,x}^{6}(J \times \mathbf{R})}^{4} \| u \|_{L_{t,x}^{6}(J \times \mathbf{R})} \| Iu \|_{L_{t,x}^{6}(J \times \mathbf{R})} \| Iu \|_{L_{t,x}^{\infty}(J \times \mathbf{R})}^{4}
\lesssim \frac{1}{N^{4}},
\endaligned
\end{equation}

\noindent and the proof of theorem $\ref{t3.1}$ is complete. $\Box$

\section{Morawetz estimates}
\begin{theorem}\label{t4.1}
Let $u$ be the solution to the nonlinear Schr{\"o}dinger equation in one dimension,

\begin{equation}\label{4.1}
i u_{t} + \Delta u = |u|^{4} u.
\end{equation}

\noindent Then

\begin{equation}\label{4.2}
\| Iu \|_{L_{t,x}^{8}([0, T] \times \mathbf{R})}^{8} \lesssim \| Iu \|_{L_{t}^{\infty} \dot{H}_{x}^{1}([0, T] \times \mathbf{R})} \| u_{0} \|_{L^{2}(\mathbf{R})}^{7} + \sum_{J_{k}} \frac{1}{N^{2-}} \| \langle \nabla \rangle Iu \|_{S^{0}(J_{k} \times \mathbf{R})}^{12},
\end{equation}

\noindent where $[0, T] = \cup_{k} J_{k}$.
\end{theorem}

\noindent \emph{Proof:} We start with the case $I = 1$. We will use the method found in \cite{DPST1} and \cite{CHVZ}. Let $\omega(t,z) : \mathbf{R} \times \mathbf{R}^{4} \rightarrow \mathbf{C}$,

\begin{equation}\label{4.3}
\omega(t,z) = u(t,x_{1}) u(t,x_{2}) u(t,x_{3}) u(t, x_{4}),
\end{equation}

\noindent where $u(t,x)$ is a solution to $(\ref{1.1})$. Then $\omega(t,z)$ obeys the equation

\begin{equation}\label{4.4}
i \omega_{t} + \Delta_{z} \omega = (\sum_{i = 1}^{4} |u(t,x_{i})|^{4}) \omega(t,z) = \mathcal N.
\end{equation}

\noindent Next, define the interaction Morawetz quantity,

\begin{equation}\label{4.5}
M_{a}(t) = 2 \int_{\mathbf{R}^{4}} \partial_{j} a(z) Im(\overline{\omega(t,z)} \partial_{j} \omega(t,z)) dz,
\end{equation}

\noindent following the convention that repeated indices are summed. Let

\begin{equation}\label{4.6}
T_{0j}(t,z) = 2 Im (\overline{\omega(t,z)} \partial_{j} \omega(t,z)),
\end{equation}

\begin{equation}\label{4.7}
L_{jk}(t,z) = -\partial_{jk} (|\omega(t,z)|^{2}) + 4 Re(\overline{\partial_{j} \omega} \partial_{k} \omega)(t,z),
\end{equation}

\begin{equation}\label{4.8}
\partial_{t} T_{0j} + \partial_{k} L_{jk} = 2 \{ \mathcal N, \omega \}_{p}^{j}.
\end{equation}

\begin{equation}\label{4.9}
\int_{0}^{T} \int_{\mathbf{R}^{4}} \partial_{t} \partial_{j} a(z) Im(\overline{\omega(t,z)} \partial_{j} \omega(t,z)) dz dt
\end{equation}

\begin{equation}\label{4.10}
= \int_{0}^{T} \int_{\mathbf{R}^{4}} \partial_{j} a(z) \partial_{jkk} (|\omega(t,z)|^{2}) dz dt
\end{equation}

\begin{equation}\label{4.11}
- 4 \int_{0}^{T} \int_{\mathbf{R}^{4}} \partial_{j} a(z) \partial_{k} Re(\overline{\partial_{j} \omega} \partial_{k} \omega)(t,z) dz dt
\end{equation}

\begin{equation}\label{4.12}
+ 2 \int_{0}^{T} \int a_{j}(z) \{ \mathcal N, \omega \}_{p}^{j} dz dt.
\end{equation}

\noindent Now, evaluate each term separately. Make a change of variables, $y = Az$, where

\begin{equation}\label{4.13}
A = \frac{1}{2} \left(
                  \begin{array}{cccc}
                    1 & 1 & 1 & 1 \\
                    1 & 1 & -1 & -1 \\
                    1 & -1 & 1 & -1 \\
                    -1 & 1 & 1 & -1 \\
                  \end{array}
                \right)
\end{equation}

\noindent is an orthonormal matrix with inverse

\begin{equation}\label{4.14}
A^{-1} = \frac{1}{2} \left(
                       \begin{array}{cccc}
                         1 & 1 & 1 & -1 \\
                         1 & 1 & -1 & 1 \\
                         1 & -1 & 1 & 1 \\
                         1 & -1 & -1 & -1 \\
                       \end{array}
                     \right).
\end{equation}

\noindent In the new variables, let $a(y) = (y_{2}^{2} + y_{3}^{2} + y_{4}^{2})^{1/2}$, $$-\Delta \Delta a(y) = 4 \pi \delta(y_{2}, y_{3}, y_{4}).$$

\begin{equation}\label{4.15}
A^{-1} \left(
         \begin{array}{c}
           y_{1} \\
           0 \\
           0 \\
           0 \\
         \end{array}
       \right) = \frac{1}{2} \left(
                               \begin{array}{c}
                                 y_{1} \\
                                 y_{1} \\
                                 y_{1} \\
                                 y_{1} \\
                               \end{array}
                             \right).
\end{equation}

\noindent Therefore, integrating $(\ref{4.10})$ by parts,

\begin{equation}\label{4.16}
(\ref{4.10}) = \int_{0}^{T} \int (-\Delta \Delta a(y)) |\omega(t,z)|^{2} dz dt =  8 \pi \int_{0}^{T} \int_{\mathbf{R}} |u(t,x)|^{8} dx dt.
\end{equation}

\noindent $\partial_{jk} a(z)$ is a positive semidefinite matrix, so integrating $(\ref{4.11})$ by parts,

\begin{equation}\label{4.17}
4 \int_{0}^{T} \int_{\mathbf{R}^{4}} (\partial_{jk} a(z)) Re (\overline{\partial_{j} \omega} \partial_{k} \omega)(t,z) \geq 0.
\end{equation}

\noindent Finally, for $(\ref{4.12})$,

\begin{equation}\label{4.18}
\{ \mathcal N, \omega \}_{p}^{j}(t,z) = -2 |\omega(t,z)|^{2} \partial_{j}(\sum_{i = 1}^{4} |u(t,x_{i})|^{4}).
\end{equation}

$$|\omega(t,z)|^{2} \partial_{j}(\sum_{i = 1}^{4} |u(t,x_{i})|^{4}) = \frac{2}{3} \sum_{j = 1}^{4} \partial_{j} (|u(t,x_{j})|^{6} (|\prod_{i \neq j} |u(t,x_{i})|^{2})),$$

\noindent so integrating $(\ref{4.12})$ by parts,

\begin{equation}\label{4.19}
-2 \int_{0}^{T} \int a_{j}(z) \{ \mathcal N, \omega \}_{p}^{j}(t,z) dz dt = \frac{4}{3} \int_{0}^{T} \int a_{jj}(z) |\omega(t,z)|^{2} |u(t,x_{j})|^{4} dz dt \geq 0.
\end{equation}

\noindent Therefore,

\begin{equation}\label{4.20}
\aligned
\int_{0}^{T} \int_{\mathbf{R}} |u(t,x)|^{8} dx dt \lesssim |\int_{0}^{T} \int \partial_{t} (a_{j}(z) Im[\overline{\omega(t,z)} \partial_{j} \omega(t,z)] dz dt| \\ \lesssim \| u \|_{L_{t}^{\infty} \dot{H}^{1/2}([0, T] \times \mathbf{R})}^{2} \| u_{0} \|_{L^{2}(\mathbf{R})}^{6}.
\endaligned
\end{equation}

\noindent Next, we prove an almost Morawetz estimate (see \cite{CGT}, \cite{CR}, \cite{D} for the two dimensional case; \cite{CGT}, \cite{DPST}, \cite{CHVZ} for discussion of the one dimensional case).\vspace{5mm}

\noindent If $u(t,x)$ solves $(\ref{4.1})$, then $Iu(t,x)$ solves

\begin{equation}\label{4.21}
i Iu_{t}(t,x) + \Delta Iu(t,x) = I(|u(t,x)|^{4} u(t,x)) = \mathcal N.
\end{equation}

\noindent Split the nonlinearity into "good" and "bad" pieces, $\mathcal N = \mathcal N_{g} + \mathcal N_{b}$.

\begin{equation}\label{4.22}
\mathcal N_{g} = \sum_{i = 1}^{4} |Iu(t,x_{i})|^{4} Iu(t,x_{i}) \prod_{j \neq i} Iu(t,x_{j}),
\end{equation}

\begin{equation}\label{4.23}
\mathcal N_{b} = \sum_{i = 1}^{4} [I(|u(t,x_{i})|^{4} u(t, x_{i})) - |Iu(t,x_{i})|^{4} (Iu(t,x_{i}))] \prod_{j \neq i} Iu(t, x_{j}).
\end{equation}

\noindent Let $\omega(t,z) = Iu(t,x_{1}) Iu(t,x_{2}) Iu(t,x_{3}) Iu(t,x_{4})$, performing the same analysis will split $$\int_{0}^{T} \int_{\mathbf{R}^{4}} \partial_{t} a_{j}(z) T_{0j}(t,z) dz dt,$$ into a sum of terms of the form $(\ref{4.10})$, $(\ref{4.11})$, and $(\ref{4.12})$. If $\mathcal N = \mathcal N_{g}$, then the previous analysis would carry over identically. Indeed,

\begin{equation}\label{4.24}
\int_{0}^{T} \int_{\mathbf{R}^{4}} (-\Delta \Delta a(z)) |\omega(t,z)|^{2} dz dt = 8 \pi \int_{0}^{T} \int |Iu(t,x)|^{8} dx dt.
\end{equation}

\begin{equation}\label{4.25}
4 \int_{0}^{T} \int_{\mathbf{R}^{4}} (\partial_{jk} a(z)) Re(\overline{\partial_{j} \omega} \partial_{k} \omega)(t,z) dz dt \geq 0.
\end{equation}

\begin{equation}\label{4.26}
2 \int_{0}^{T} \int a_{j}(z) \{ \mathcal N_{g}, \omega \}_{p}^{j}(t,z) dz dt = \frac{4}{3} \int_{0}^{T} \int a_{jj}(z) |\omega(t,z)|^{2} |Iu(t,x_{j})|^{4} dz \geq 0.
\end{equation}

\noindent Therefore,

\begin{equation}\label{4.27}
\int_{0}^{T} \int |u(t,x)|^{8} dx dt \lesssim \| Iu \|_{L_{t}^{\infty} \dot{H}_{x}^{1}([0, T] \times \mathbf{R})} \| u_{0} \|_{L^{2}(\mathbf{R})}^{7} + |\int_{0}^{T} \int a_{j}(z) \{ \mathcal N_{b}, \omega \}_{p}^{j}(t,z) dz dt|.
\end{equation}

\noindent To analyze the remainder $\mathcal N_{b}$, first consider a term of the form

\begin{equation}\label{4.28}
\int_{J_{k}} \int a_{j}(z) \overline{\mathcal N_{b}(t,z)} \partial_{j} \omega(t,z) dz dt.
\end{equation}

\noindent Recall $(\ref{4.23})$, without loss of generality let $i = 1$ and estimate

\begin{equation}\label{4.28.1}
\aligned
\int_{J_{k}} \int &a_{j}(z) \overline{[I(|u(t,x_{1})|^{4} u(t,x_{1})) - |Iu(t,x_{1})|^{4} Iu(t,x_{1})] Iu(t,x_{2}) Iu(t,x_{3}) Iu(t,x_{4})} \\ &\times \partial_{j}(Iu(t,x_{1}) Iu(t,x_{2}) Iu(t,x_{3}) Iu(t,x_{4})) dx dt.
\endaligned
\end{equation}

\noindent Because $a_{j}(z) \in L^{\infty}$,

\begin{equation}\label{4.29}
(\ref{4.28}) \lesssim \| I(|u(t,x)|^{4} u(t,x)) - |Iu(t,x)|^{4} Iu(t,x) \|_{L_{t}^{1} L_{x}^{2}(J_{k} \times \mathbf{R})} \| \langle \nabla \rangle Iu \|_{L_{t}^{\infty} L_{x}^{2}(\mathbf{R})}^{7}.
\end{equation}

\noindent To evaluate $$ I(|u(t,x)|^{4} u(t,x)) - |Iu(t,x)|^{4} Iu(t,x),$$ make a Littlewood - Paley partition of unity. Let

\begin{equation}\label{4.30}
F(t,\xi) = \int_{\xi = \xi_{1} + ... + \xi_{5}} [1 - \frac{m(\xi_{1} + ... + \xi_{5})}{m(\xi_{1}) ... m(\xi_{5})}] \widehat{Iu}(t,\xi_{1}) \widehat{\overline{Iu}}(t,\xi_{2}) \widehat{Iu}(t,\xi_{3}) \widehat{\overline{Iu}}(t,\xi_{4}) \widehat{Iu}(t,\xi_{5}).
\end{equation}

\noindent Without loss of generality, let $N_{1} \geq N_{2} \geq N_{3} \geq N_{4} \geq N_{5}$. Consider several cases separately.\vspace{5mm}

\noindent \textbf{Case 1: $N_{1} << N$} In this case, $$|1 - \frac{m(\xi_{1} + ... + \xi_{5})}{m(\xi_{1}) ... m(\xi_{5})}| \equiv 0.$$

\noindent \textbf{Case 2: $N_{1} \gtrsim N >> N_{2}$} In this case, by the fundamental theorem of calculus, $$|1 - \frac{m(\xi_{1} + ... + \xi_{5})}{m(\xi_{1}) ... m(\xi_{5})}| \lesssim \frac{N_{2}}{N_{1}}.$$

$$(\ref{4.28}) \lesssim \sum_{N \lesssim N_{1}} \sum_{N_{5} \leq N_{4} \leq N_{3} \leq N_{2} << N} \frac{1}{N_{1}^{2}} \| P_{N_{1}} \langle \nabla \rangle Iu \|_{L_{t}^{\infty} L_{x}^{2}(J_{k} \times \mathbf{R})} \| P_{N_{2}} \langle \nabla \rangle Iu \|_{L_{t}^{4} L_{x}^{\infty}(J_{k} \times \mathbf{R})}$$

$$\times  \| P_{N_{3}} \langle \nabla \rangle Iu \|_{L_{t}^{4} L_{x}^{\infty}(J_{k} \times \mathbf{R})}  \| P_{N_{4}} \langle \nabla \rangle Iu \|_{L_{t}^{4} L_{x}^{\infty}(J_{k} \times \mathbf{R})}  \| P_{N_{5}} \langle \nabla \rangle Iu \|_{L_{t}^{4} L_{x}^{\infty}(J_{k} \times \mathbf{R})}$$

$$\lesssim \frac{1}{N^{2-}} \| \langle \nabla \rangle Iu \|_{S^{0}(J \times \mathbf{R})}^{5}.$$

\noindent \textbf{Case 3: $N_{1} \geq N_{2} \gtrsim N$} In this case, crudely estimate $$|1 - \frac{m(\xi_{1} + ... + \xi_{5})}{m(\xi_{1}) ... m(\xi_{5})}| \lesssim \frac{1}{m(N_{1}) m(N_{2}) m(N_{3}) m(N_{4}) m(N_{5})}.$$

$$(\ref{4.28}) \lesssim \sum_{N \lesssim N_{2} \leq N_{1}} \| P_{N_{1}} u \|_{L_{t}^{\infty} L_{x}^{2}(J_{k} \times \mathbf{R})} \| P_{N_{2}} u \|_{L_{t}^{4} L_{x}^{\infty}(J_{k} \times \mathbf{R})}$$

$$\times \sum_{N_{5} \leq N_{4} \leq N_{3} \leq N_{2}} \| P_{N_{3}} u \|_{L_{t}^{4} L_{x}^{\infty}(J_{k} \times \mathbf{R})}  \| P_{N_{4}} u \|_{L_{t}^{4} L_{x}^{\infty}(J_{k} \times \mathbf{R})}  \| P_{N_{5}} u \|_{L_{t}^{4} L_{x}^{\infty}(J_{k} \times \mathbf{R})}$$

$$\lesssim \sum_{N \lesssim N_{2} \leq N_{1}} \frac{1}{N^{1 - s} N_{1}^{s}} \frac{1}{N^{1 - s} N_{2}^{s}}$$ $$\times \sum_{N_{5} \leq N_{4} \leq N_{3}} \frac{1}{m(N_{3}) N_{3}} \frac{1}{m(N_{4}) N_{4}} \frac{1}{m(N_{5}) N_{5}} \| \langle \nabla \rangle Iu \|_{S^{0}(J_{k} \times \mathbf{R})}^{5}$$

$$\lesssim \frac{1}{N^{2-}} \| \langle \nabla \rangle Iu \|_{S^{0}(J_{k} \times \mathbf{R})}^{5}.$$

\noindent For a term of the form

\begin{equation}\label{4.31}
\int_{J_{k}} \int_{\mathbf{R}^{4}} a_{j}(z) (\partial_{j} \overline{\mathcal N_{b}(t,z)}) \omega(t,z) dz dt,
\end{equation}

\noindent integrating by parts rewrites this term as a term of the form $(\ref{4.28})$ plus a term of the form

\begin{equation}\label{4.32}
\int_{J_{k}} \int_{\mathbf{R}^{4}} a_{jj}(z) (\overline{\mathcal N_{b}(t,z)}) \omega(t,z) dz dt.
\end{equation}

\noindent By $(\ref{4.23})$ and the estimates on $(\ref{4.30})$,

\begin{equation}\label{4.33}
\| \mathcal N_{b}(t,z) \|_{L_{t}^{1} L_{x}^{2}(J_{k} \times \mathbf{R}^{4})} \lesssim \frac{1}{N^{2-}} \| \langle \nabla \rangle Iu \|_{S^{0}(J_{k} \times \mathbf{R})}^{8}.
\end{equation}

\noindent $\Delta a$ does not lie in $L^{\infty}(\mathbf{R}^{4})$, rather,

\begin{equation}\label{4.34}
a_{jj}(y) \lesssim \frac{1}{(y_{2}^{2} + y_{3}^{2} + y_{4}^{2})^{1/2}}.
 \end{equation}

 \noindent Integrating $$\frac{1}{(y_{2}^{2} + y_{3}^{2} + y_{4}^{2})^{1/2}} Iu(t,x_{1}) Iu(t,x_{2}) Iu(t,x_{3}) Iu(t,x_{4})$$ along the tube $y_{2}^{2} + y_{3}^{2} + y_{4}^{2} \leq 1$, using $(\ref{4.15})$,

$$\int_{-\infty}^{\infty} \int_{y_{2}^{2} + y_{3}^{2} + y_{4}^{2} \leq 1}  \frac{1}{y_{2}^{2} + y_{3}^{2} + y_{4}^{2}} |Iu(t, \frac{y_{1} + y_{2} + y_{3} - y_{4}}{2})|^{2} |Iu(t, \frac{y_{1} + y_{2} - y_{3} + y_{4}}{2})|^{2}$$

$$\times  |Iu(t, \frac{y_{1} - y_{2} + y_{3} + y_{4}}{2})|^{2} |Iu(t, \frac{y_{1} - y_{2} - y_{3} - y_{4}}{2})|^{2} dy_{1} dy_{2} dy_{3} dy_{4}$$

$$\lesssim (\int_{0}^{1} \frac{r^{2}}{r^{2}} dr) \| Iu \|_{L_{t}^{\infty} L_{x}^{8}(J_{k} \times \mathbf{R})}^{8} \lesssim \| \langle \nabla \rangle Iu \|_{S^{0}(J_{k} \times \mathbf{R})}^{8},$$

\noindent so $$\| a_{jj}(y) \omega(t,y) \|_{L_{t}^{\infty} L_{y}^{2}(J_{k} \times \mathbf{R}^{4})} \lesssim \| \langle \nabla \rangle Iu \|_{S^{0}(J_{k} \times \mathbf{R})}^{4}.$$

\noindent On the other hand, when $y_{2}^{2} + y_{3}^{2} + y_{4}^{2} \geq 1$, $\partial_{jj} a(z)$ is bounded and

\begin{equation}\label{4.35}
\| Iu(t,x_{1}) Iu(t,x_{2}) Iu(t,x_{3}) Iu(t,x_{4}) \|_{L_{t}^{1} L_{x}^{2}(J_{k} \times \mathbf{R}^{4})} \lesssim \| \langle \nabla \rangle Iu \|_{S^{0}(J_{k} \times \mathbf{R})}^{4}.
\end{equation}

\noindent Since $$\| I(|u(t,x_{1})|^{4} u(t,x_{1})) - |Iu(t,x_{1})|^{4} Iu(t,x_{1}) \|_{L_{t}^{1} L_{x}^{2}(J \times \mathbf{R})} \lesssim \frac{1}{N^{2-}} \| \langle \nabla \rangle Iu \|_{S^{0}(J_{k} \times \mathbf{R})}^{5},$$ theorem $\ref{t4.1}$ is proved.

\section{Proof of Theorem $\ref{t1.1}$}

\begin{theorem}\label{t5.1}
$(\ref{1.1})$ is globally well-posed for $u_{0} \in H^{s}(\mathbf{R})$, $s > \frac{8}{29}$.
 \end{theorem}

 \noindent \emph{Proof:}

\begin{equation}\label{5.1}
\aligned
\int |\nabla Iu_{0}(x)|^{2} dx \lesssim N^{2(1 - s)} \| u_{0} \|_{H^{s}(\mathbf{R})}. \\
\int |Iu_{0}(x)|^{6} dx \lesssim N^{2 - 6s} \| u_{0} \|_{H^{s}(\mathbf{R})}^{6}.
\endaligned
\end{equation}

\noindent If $u(t,x)$ solves $(\ref{1.1})$ on $[0, T_{0}]$, then rescaling,

\begin{equation}\label{5.2}
\frac{1}{\lambda^{1/2}} u(\frac{t}{\lambda^{2}}, \frac{x}{\lambda})
\end{equation}

\noindent solves $(\ref{1.1})$ on $[0, \lambda^{2} T_{0}]$, we will call the rescaled solution $u_{\lambda}(t,x)$. Choose $\lambda \sim N^{(1 - s)/s}$ so that $E(Iu_{\lambda}(0)) = 1/2$. Let

\begin{equation}\label{5.3}
W = \{ t : E(Iu_{\lambda}(t)) \leq \frac{9}{10} \}.
\end{equation}

\noindent W is closed by the dominated convergence theorem and nonempty since $0 \in W$. To prove $W = [0, \lambda^{2} T_{0}]$, it suffices to prove $W$ is open in $[0, \lambda^{2} T_{0}]$. Suppose $W = [0, T]$, then by continuity there exists $\delta > 0$ such that $E(Iu_{\lambda}(t)) \leq 1$ on $[0, T + \delta]$.

\begin{lemma}\label{l5.1}
\noindent When $s > 1/4$,

\begin{equation}\label{5.4}
\| Iu_{\lambda} \|_{L_{t,x}^{8}([0, T + \delta] \times \mathbf{R})}^{8} \leq \frac{3}{2} C_{0},
\end{equation}

\noindent for some $C_{0}(T_{0})$.
\end{lemma}

\noindent \emph{Proof:} Let $\tau = \sup \{ T_{\ast} \in [0, T + \delta] : \| Iu_{\lambda} \|_{L_{t,x}^{8}([0, T_{\ast}] \times \mathbf{R})}^{8} \leq \frac{3}{2} C_{0} \}$. If $\tau < T + \delta$, there exists $\delta' > 0$ such that

\begin{equation}\label{5.5}
\| Iu_{\lambda} \|_{L_{t,x}^{8}([0, \tau + \delta'] \times \mathbf{R})}^{8} \leq 2 C_{0}.
\end{equation}

\noindent Recall $E(Iu_{\lambda}(t)) \leq 1$ on $[0, \tau + \delta']$,

\begin{equation}\label{5.6}
\| Iu_{\lambda} \|_{L_{t}^{16/3} L_{x}^{8}([0, \tau + \delta'] \times \mathbf{R})} \leq (\lambda^{2} T_{0})^{1/16} \| Iu_{\lambda} \|_{L_{t,x}^{8}([0, \tau + \delta'] \times \mathbf{R})} \leq (\lambda^{2} T_{0})^{1/16} (2C_{0})^{1/8}.
\end{equation}

\noindent Partition $[0, \tau + \delta']$ into $$\frac{(\lambda^{2} T_{0})^{1/3} (2C_{0})^{2/3}}{\epsilon^{16/3}}$$ subintervals such that

\begin{equation}\label{5.7}
\| Iu_{\lambda} \|_{L_{t}^{16/3} L_{x}^{8}([0, \tau + \delta'] \times \mathbf{R})} \leq \epsilon
\end{equation}

\noindent on each subinterval. Then apply the almost Morawetz estimate,

\begin{equation}\label{5.8}
\aligned
\| Iu_{\lambda} \|_{L_{t,x}^{8}([0, \tau + \delta'] \times \mathbf{R})}^{8} \leq C \| u_{0} \|_{L^{2}(\mathbf{R})}^{7} \| Iu_{\lambda} \|_{L_{t}^{\infty} \dot{H}^{1}([0, \tau + \delta'] \times \mathbf{R})} + \sum_{k} \frac{\| \langle \nabla \rangle Iu \|_{S^{0}(J_{k} \times \mathbf{R})}^{12}}{N^{2-}} \\
\leq C_{0} + C \frac{(\lambda^{2} T_{0})^{1/3} (2C_{0})^{2/3}}{\epsilon^{16/3} N^{2-}}
\leq  C_{0} + C \frac{N^{\frac{2(1 - s)}{3s}} T_{0}^{1/3} (2C_{0})^{2/3}}{\epsilon^{16/3} N^{2-}}
\leq \frac{3}{2} C_{0},
\endaligned
\end{equation}

\noindent when $N$ is sufficiently large, as long as $\frac{2}{3} (\frac{1 - s}{s}) < 2$, or $s > 1/4$. This proves the lemma. $\Box$\vspace{5mm}

\noindent Returning to the theorem,

\begin{equation}\label{5.9}
\aligned
\| Iu_{\lambda} \|_{L_{t,x}^{8}([0, T + \delta] \times \mathbf{R})}^{8} \leq \frac{3}{2} C_{0}, \\
\| Iu_{\lambda} \|_{L_{t}^{16/3} L_{x}^{8}([0, T + \delta] \times \mathbf{R})} \lesssim \lambda^{1/8} T_{0}^{1/16}.
\endaligned
\end{equation}

\noindent Partition $[0, T + \delta]$ into $\lesssim (\lambda^{2} T_{0})^{1/3}$ subintervals. We will call these the little intervals. Take the union of the first $N^{1/2-}$ little subintervals, and call this big interval $$J_{1} = \cup_{l = 1}^{N^{1/2-}} J_{1,l}.$$ Take the union of the next $N^{1/2-}$ subintervals and call this big interval $J_{2}$, and so on.

\begin{equation}\label{5.10}
[0, T + \delta] = \bigcup_{k = 1}^{\frac{(\lambda^{2} T_{0})^{1/3}}{N^{1/2-}}} J_{k} = \bigcup_{k = 1}^{\frac{(\lambda^{2} T_{0})^{1/3}}{N^{1/2-}}} \bigcup_{l = 1}^{N^{1/2-}} J_{k,l}.
\end{equation}

\begin{lemma}\label{l5.2}
\noindent Suppose $E(Iu(t)) \leq 1$ on $J_{k}$.

\begin{equation}\label{5.11}
\sup_{t_{1}, t_{2} \in J_{k}} |E(Iu(t_{1})) - E(Iu(t_{2})) | \lesssim \frac{1}{N^{5/4-}}.
\end{equation}
\end{lemma}

\noindent \emph{Proof:} By theorem $(\ref{t3.1})$,

\begin{equation}\label{5.12}
\sup_{t_{1}, t_{2} \in J_{k,l}} |E(Iu(t_{1})) - E(Iu(t_{2}))| \lesssim \frac{1}{N^{3/2-}} \| P_{> cN} \nabla Iu \|_{L_{t}^{4} L_{x}^{\infty}(J_{k,l} \times \mathbf{R})}^{2} + \frac{1}{N^{2-}}
\end{equation}

$$\sum_{l = 1}^{N^{1/2}} \| P_{> cN} \nabla Iu \|_{L_{t}^{4} L_{x}^{\infty}(J_{k,l} \times \mathbf{R})}^{2} \lesssim N^{1/4} \| P_{> cN} \nabla Iu \|_{L_{t}^{4} L_{x}^{\infty}(J_{k} \times \mathbf{R})}^{2}.$$

\noindent Let $J_{k} = [a, b]$, $J_{k,m} = [a_{m}, b_{m}]$, $a_{m + 1} = b_{m}$, $a_{1} = a$, $b_{N^{1/2-}} = b$. On each little subinterval, perform the linear-nonlinear decomposition in theorem $\ref{t2.1}$. The solution on $[a_{m}, b_{m}]$ is of the form $$e^{i(t - a_{m}) \Delta} u(a_{m}) + u_{m}^{nl}(t).$$ By induction,

\begin{equation}\label{5.12.1}
e^{i(t - a_{m}) \Delta} u(a_{m}) = e^{it \Delta} u_{0} + \sum_{j = 1}^{m - 1} e^{i(t - a_{j}) \Delta} u_{j}^{nl}(b_{j}).
\end{equation}

\begin{equation}\label{5.13}
\| \nabla e^{i(t - a) \Delta} Iu(a) \|_{L_{t}^{4} L_{x}^{\infty}(J_{k} \times \mathbf{R})} \lesssim 1.
\end{equation}

\begin{equation}\label{5.14}
\sum_{m = 1}^{N^{1/2-}} \| \nabla e^{i(t - a) \Delta} I u(b_{m})^{nl} \|_{L_{t}^{4} L_{x}^{\infty}(J_{k} \times \mathbf{R})} \lesssim 1.
\end{equation}

\begin{equation}\label{5.15}
\sum_{l = 1}^{N^{1/2-}} \| \nabla Iu^{nl}(t) \|_{L_{t}^{4} L_{x}^{\infty}(J_{k,l} \times \mathbf{R})}^{4} \lesssim \frac{1}{N^{3/2-}}.
\end{equation}

\noindent Therefore, $$\| \nabla P_{> cN} Iu \|_{L_{t}^{4} L_{x}^{\infty}(J_{k} \times \mathbf{R})} \lesssim 1.$$ Plugging this back in to $(\ref{5.11})$, $$\sum_{l = 1}^{N^{1/2-}} \sup_{t_{1}, t_{2} \in J_{k,l}} |E(Iu(t_{1})) - E(Iu(t_{2}))| \lesssim \frac{1}{N^{5/4-}} + \frac{1}{N^{3/2-}},$$

\noindent which proves the lemma. $\Box$\vspace{5mm}

\noindent Returning to the theorem once again, if $(\lambda^{2} T_{0})^{1/3} << N^{7/4-}$, for $t \in [0, T + \delta]$, or if $s > \frac{8}{29}$,

\begin{equation}\label{5.16}
E(Iu(t)) \leq \frac{1}{2} + C \frac{(\lambda^{2} T_{0})^{1/3}}{N^{7/4-}} = \frac{1}{2} + C \frac{N^{\frac{2(1 - s)}{3s}} T_{0}^{1/3}}{N^{7/4-}} \leq \frac{9}{10},
\end{equation}

\noindent for sufficiently large N, proving the theorem. We have to take $N \sim T_{0}^{\frac{4s}{29s - 8}}$. Since $\lambda \sim N^{\frac{1 - s}{s}}$, and $$\| u(t) \|_{H^{s}(\mathbf{R})} = \lambda^{s} \| u_{\lambda}(\frac{t}{\lambda^{2}}) \|_{H^{s}(\mathbf{R})},$$

\begin{equation}\label{5.17}
\sup_{t \in [0, T]} \| u(t) \|_{H^{s}(\mathbf{R})} \lesssim (1 + T)^{\frac{4(1 - s)s}{29s - 8}+}.
\end{equation}

 $\Box$

\newpage
\nocite*
\bibliographystyle{plain}
\bibliography{n=1}

\end{document}